\documentclass[11pt]{amsart}

\usepackage{amsmath}
\usepackage{amssymb}
\usepackage{mathrsfs}

\usepackage{graphicx}

\newtheorem{lemat}{Lemma}
\newtheorem{theorem}{Theorem}
\newtheorem{fact}{Fact}
\newtheorem{definition}{Definition}

\newcommand{\nat}{\mathbb{N}}
\newcommand{\set}[2]{\{#1,\ldots,#2\}}
\newcommand{\nt}{\mathcal{T}}
\newcommand{\nn}{\mathcal{N}}

\newcommand{\niceimage}[3]{\begin{figure}[h]\begin{center}\includegraphics[width=#3\textwidth]{#1}\caption{#2}\end{center}\end{figure}}

\begin{document}

\author[Michal Adamaszek]{Micha{\l} Adamaszek$^*$}
\thanks{$^*$ Warsaw University, email: aszek@mimuw.edu.pl}
\title{Efficient enumeration of graceful permutations}

\begin{abstract}
A graceful $n$-permutation is a graceful labeling of an $n$-vertex path $P_n$. In this paper we improve the asymptotic lower bound on the number of such permutations from $\Omega((5/3)^n)$ to $\Omega(2.37^n)$. This is a computer-assisted proof based on an effective algorithm that enumerates graceful $n$-permutations. Our algorithm is also presented in detail.
\end{abstract}


\maketitle

\section{Graceful graphs and permutations}

Let $G=(V,E)$ be an undirected graph with $|V| = n$ and $|E| = m$. We say that a vertex labeling $f:V\mapsto \nat$ together with an edge labeling $g:E\mapsto \nat$ are a \emph{graceful labeling} of $G$ if:
\begin{itemize}
	\item $f(V) \subset \set{0}{m}$ and $f$ is one-to-one (injective)
	\item $g(E) = \set{1}{m}$
	\item $g(uv) = |f(u)-f(v)|$ for every two vertices $u,v\in V$ such, that $uv\in E$
\end{itemize}

Graceful labelings of graphs have received a lot of attention; see \cite{Gal05} for an extensive survey. In this paper we concentrate on the single case when $G = P_n$ is the $n$-vertex path. Note, that in this case $m=n-1$, thus the vertex labels are in bijection with the set $\set{0}{n-1}$. This justifies the following definition:
\begin{definition}A permutation $[\sigma(0),\ldots,\sigma(n-1)]$ of the set $\{0,1,\ldots,n-1\}$ is a \emph{graceful $n$-permutation} if
	$$\{|\sigma(1)-\sigma(0)|,|\sigma(2)-\sigma(1)|,\ldots,|\sigma(n-1)-\sigma(n-2)|\}=\set{1}{n-1}$$
\end{definition}
For instance, $[0,6,1,5,2,4,3]$ is a graceful $7$-permutation. The values of a graceful $n$-permutation can be identified with the vertex labels in some graceful labeling of $P_n$ and vice versa. We shall use these notions interchangeably.

Denote by $G(n)$ the number of graceful $n$-permutations. The sequence $G(n)$ is not well known, not even asymptotically. It has number A006967 in the Sloane's On-line Encyclopedia of Integer Sequences (\cite{OEIS}) where the first 20 terms are listed. Its growth is exponential as shown in \cite{Klo95} and \cite{AlS03}. In the latter the best known estimate, $G(n)=\Omega((\frac{5}{3})^n)$ is proved. Here we shall improve this result by proving the following:
\begin{theorem}\label{theo1} $G(n)=\Omega(2.37^n)$\end{theorem}

This paper is organized as follows. In the next section we introduce a recursive algorithm for the computation of $G(n)$. Next we observe how its efficiency can be vastly improved using some knowledge of the structure of graceful permutations. In section 4 we use the computational data to prove Theorem \ref{theo1}. Some closing remarks are included in section 5.

\section{The search tree}

We shall generate (and count) graceful $n$-permutations by the following recursive search (think of path labelings for now): the edge label $n-1$ can only appear as $|0-(n-1)|$, therefore the vertices with labels $0$ and $n-1$ must be neighbours. Moving along, the next free edge label $n-2$ can be induced as $|0-(n-2)|$ or $|1-(n-1)|$, so either $0$ and $n-2$ or $1$ and $n-1$ must be connected. This obvious procedure continues with further edge labels down to $1$. Of course we can only test adding a certain edge if it does not conflict with the path structure of the created graph, i.e. if what has been constructed so far is a collection of paths. 

Figure 1 shows half of the \emph{search tree $\nt_n$} obtained for $n=7$. The nodes of the tree will be referred to as \emph{partial permutations}. The \emph{level}, indicated in the left column, is the edge label just added. The vertex labels inducing that edge label are underlined in each node. The $8$ nodes on level $1$ give rise to $16$ graceful $7$-permutations (each can be read in the given order or backwards), therefore $G(7)=32$ (because the other half of the tree looks just the same).

\niceimage{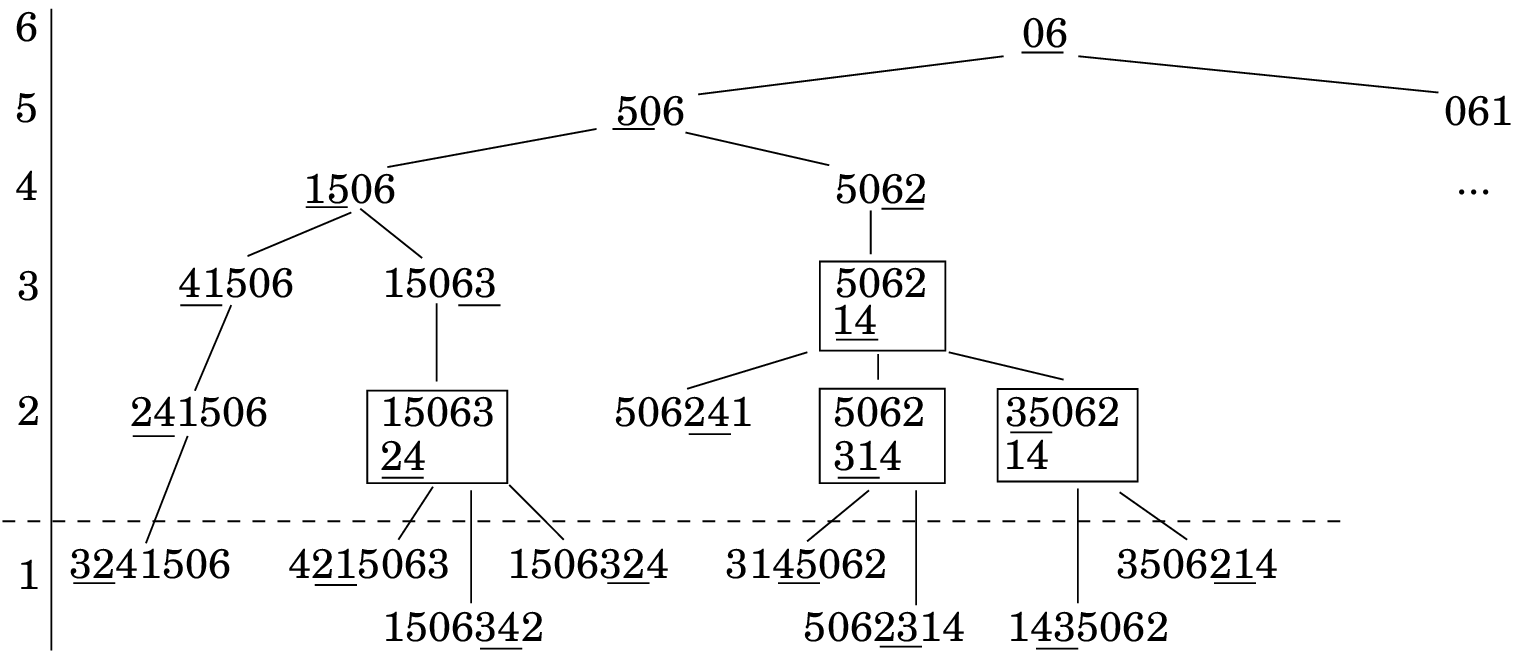}{Half of the search tree $\nt_7$}{0.9}

Now we describe how the partial permutations are represented, so that the expansion from a node on level $k$ to the nodes on level $k-1$ can be done quickly. Each node holds two arrays:

\begin{itemize}
	\item \verb+free[0..n-1]+ --- for each $u$ the number \verb+free[u]+ is the number of ,,free slots'' of a vertex labeled $u$. Initially \verb+free[u]:=2+ for all $u$, and the value drops down by one each time $u$ is chosen an endpoint of some edge. A vertex label with \verb+free[u]=0+ must no longer be used,
	\item \verb+forb[0..n-1]+ --- (forbidden): for every vertex label $u$, which is an endpoint of some partial path in the partial permutation \verb+forb[u]+ is the label of the other end of this path. These two vertices cannot be connected by an edge, since a cycle would appear. Initially \verb+forb[u]:=u+.
\end{itemize}

Note that with this representation we do not know the actual order of the labels in the permutation. However, they fall into three classes: yet unused (with \verb+forb[u]=2+), endpoints of partial paths (\verb+forb[u]=1+) and no longer available labels inside the paths (\verb+free[u]=0+). As for the endpoints their pairing is completely described by \verb+forb+. Note that \verb+forb[forb[u]]=u+ at all times.

Expansion is now easy: a new edge can be added between two labels iff they both have at least one free slot and are not paired by \verb+forb+. To update the \verb+forb+ array after a new edge addition note, that an edge can be added in three ways: between two yet unused labels, between an unused label and a path endpoint or between two endpoints of different paths. Supposing that the labels being connected are $u$ and $v$ the following simultaneous assignment:
\begin{center}\verb+(forb[forb[u]], forb[forb[v]]) := (forb[v], forb[u])+\end{center}
is valid in each case, which is an easy check to verify.

A straightforward recursive tree search procedure that counts graceful $n$-permutations is an obvious outcome of the above considerations. Now we shall work on efficiency.

\section{Equivalence of partial permutations}

In Figure 1 half of the search tree was omitted because it resembles the first half. More precisely, if $f:V\mapsto\nat$ and $g:E\mapsto\nat$ is a graceful labeling of a graph $G=(V,E)$ with $m$ edges then the \emph{complementary labeling} given by:
	$$\overline{f}(v) = m-f(v), \overline{g}(uv)=g(uv)$$
is again graceful. The omitted half of the tree was, in this sense, complementary to the first half so it yielded equally many graceful permutations.

Now we shall generalize this, and define an equivalence relation between the nodes on one level in the search tree.

\begin{definition}
	Let ($free_1$, $forb_1$) and ($free_2$, $forb_2$) be the arrays \texttt{free} and \texttt{forb} in two partial permutations $\nn_1$ and $\nn_2$ on the same level of the search tree $\nt_n$. We say the nodes $\nn_1$ and $\nn_2$ are \emph{equivalent} if 
	either
	$$\forall_{u}\ free_1[u] = free_2[u] \textrm{ and } \forall_{u}\ (free_1[u]=1 \Rightarrow forb_1[u]=forb_2[u])$$
	or
	$$\forall_{u}\ free_1[u] = free_2[n-1-u] \textrm{ and } \forall_{u}\ (free_1[u]=1 \Rightarrow forb_1[u]=n-1-forb_2[n-1-u])$$
\end{definition}

Less formally it says that $\nn_1$ and $\nn_2$ are equivalent if they have the same \verb+forb+-pairing of endpoints and the same set of used labels, possibly after taking the complementary labeling in one of the nodes. This is an equivalence relation with the following additional property:

\begin{fact}
If $\nn_1$ and $\nn_2$ are equivalent then the number of graceful permutations they expand to (i.e. number of leaves on level $1$ in the subtrees of $\nt_n$ rooted in $\nn_1$ and $\nn_2$ respectively) are equal.
\end{fact}
\textbf{Proof.} This follows from the remark in the previous section, that only the \verb+forb+-pairing and the set of free labels influence the expansion algorithm (exact location of inside-path labels does not matter). On the other hand complementary nodes yield symmetrical (complementary) subtrees. In either case equivalent partial permutations expand to isomorphic rooted subtrees of $\nt_n$.\qed

This observation leads to a breadth-first search of the search tree. With each node we keep its \emph{multiplicity} --- the number of nodes in its equivalence class. We only keep one representative of each class. After expansion from level $k$ to $k-1$ we group the new nodes into equivalence classes again and sum up multiplicities accordingly. The final answer $G(n)$ is the sum of multiplicities of all nodes on level $1$. Note, that comparing two nodes with respect to equivalence takes $\Theta(n)$ time, thus full comparison of new nodes during the expansion from level $k$ to $k-1$ would be expensive. To speed this up a hash table was used to keep new nodes. Observe, that the choice of the hash function is not completely arbitrary --- it must not distinguish equivalent nodes. 

Here are some sample numbers to indicate the power of the optimization thus achieved: $G(40)\approx 0.2\cdot10^{18}$ is the number of nodes at level $1$ in $\nt_{40}$. Hovewer, there are less than $3\cdot10^5$ distinct equivalence classes of partial permutations at each level, therefore at most this many partial permutations must be kept in memory and expanded at a time.

\section{The main results}

To get a lower bound on $G(n)$ we follow precisely the method of \cite{Klo95} and \cite{AlS03}. First, extend the notation $G(n)$ to:
\begin{itemize}
	\item $G(n;a)$ --- the number of graceful $n$-permutations with left endpoint $a$ \newline\mbox{($\sigma(0)=a$)}, let us call them graceful $(n;a)$-permutations
	\item $G(n;a,b)$ --- the number of graceful $n$-permutations with endpoints $a,b$ \newline ($\sigma(0)=a, \sigma(n-1)=b$), let us call them graceful $(n;a,b)$-permutations
\end{itemize}

\begin{lemat}{(\cite{Klo95}, \cite{AlS03})}
	\label{lemma}
	For any numbers $r,m,j$, $j\leq m$ we have the inequality:
	$$G(r+2m;j)\geq G(2m;j,j+m) G(r;j)$$
\end{lemat}
\textbf{Proof (sketch).} First prove that a graceful $(2m;j,j+m)$-permutation is in fact \emph{bipartite} graceful --- all edges connect large (greater or equal $m$) vertex labels with small ones. Then add $r$ to all large vertex labels in a $(2m;j,j+m)$-permutation, add $m$ to all vertex labels in a $(r;j)$-permutation and glue these two by adding an edge between $j+m+r$ and $j+m$. This yields a $(r+2m;j)$-permutation.\qed

By iterating the last inequality $k$ times we get:
	$$G(r+2km;j)\geq G(2m;j,j+m)^k G(r;j)$$
Hence, for fixed $m,j$ we have an estimate ($n$ is the variable):
	$$G(n)\geq G(n;j) = \Omega(\ (G(2m;j,j+m)^{\frac{1}{2m}})\ ^n)$$
It remains to find $m,j$, that make $\gamma_{m,j}=G(2m;j,j+m)^{\frac{1}{2m}}$ possibly large. Observations show, that for a fixed $m$ $\gamma_{m,j}$ is the biggest for $j=\lfloor\frac{m}{2}\rfloor$ and that the sequence $\gamma_{m,\lfloor\frac{m}{2}\rfloor}$ is increasing. Hence it is desirable to compute $\gamma_{m,\lfloor\frac{m}{2}\rfloor}$ form as big $m$ as possible, which is equivalent to computing $G(2m;\lfloor\frac{m}{2}\rfloor,\lfloor\frac{m}{2}\rfloor+m)$. The results so far were:
\begin{itemize}
	\item in \cite{Klo95}: $G(20;5,15)=4\,382$, hence $G(n)=\Omega(1.521^n)$
	\item in \cite{AlS03}: $G(26;6,19)=636\,408$, hence $G(n)=\Omega(1.671^n)=\Omega((\frac{5}{3})^n)$
\end{itemize}

With slight easy modifications the algorithm described in the previous section can be used to enumerate also graceful $(n;a,b)$-permutations. It was efficient enough to compute:

$$G(64;16,48) = 1\,172\,380\,428\,523\,169\,632\,220\,649$$

which in turn yields $\gamma_{32,16} = G(64;16,48)^{1/64} > (10^{24})^{1/64} > 2.37$. Eventually we get:

$$G(n)\geq G(n; 16) = \Omega(2.37^n)$$

This completes the proof of Theorem \ref{theo1}.

\section{Closing remarks}
Additionally the values of $G(n)$ have been computed for $n\leq 40$ (they have been submitted to \cite{OEIS}). The quotients $G(n+1)/G(n)$ tend to gather between $3$ and $4.5$, suggesting that the lower bound $2.37^n$ is poor. This is no surprise, because we have in fact estimated the size of only a small part of all graceful $n$-permutations, namely the bipartite graceful $(n;16)$-permutations. It also remains an open question to find an exponential upper bound on $G(n)$.

\end{document}